\documentclass{article}
\usepackage{amsmath}
\usepackage{mathrsfs}
\usepackage{txfonts}
\usepackage{bm}
\usepackage{latexsym,amsfonts,amssymb}
\begin{document}
\setcounter{page}{1}
\newtheorem{thm}{Theorem}[section]
\newtheorem{fthm}[thm]{Fundamental Theorem}
\newtheorem{dfn}[thm]{Definition}
\newtheorem{rem}[thm]{Remark}
\newtheorem{lem}[thm]{Lemma}
\newtheorem{cor}[thm]{Corollary}
\newtheorem{exa}[thm]{Example}
\newtheorem{pro}[thm]{Proposition}
\newtheorem{prob}[thm]{Problem}

\newtheorem{theorem}{Theorem}[section]
\newtheorem{definition}[theorem]{Definition}
\newtheorem{remark}[theorem]{Remark}
\newtheorem{lemma}[theorem]{Lemma}
\newtheorem{corollary}[theorem]{Corollary}
\newtheorem{example}[theorem]{Example}
\newtheorem{proposition}[theorem]{Proposition}
\newtheorem{con}[thm]{Conjecture}
\newtheorem{ob}[thm]{Observation}

\renewcommand{\theequation}{\thesection.\arabic{equation}}
\renewcommand{\thefootnote}{\fnsymbol{footnote}}

\renewcommand{\thefootnote}{\fnsymbol{footnote}}
\newcommand{\qed}{\hfill\Box\medskip}
\newcommand{\proof}{{\it Proof.\quad}}

\newcommand{\rmnum}[1]{\romannumeral #1}
\renewcommand{\abovewithdelims}[2]{%
\genfrac{[}{]}{0pt}{}{#1}{#2}}

\renewcommand{\thefootnote}{\fnsymbol{footnote}}

\title{\bf Generalized bilinear forms graphs and MDR codes\\ over residue class rings\footnote{Project 11371072
supported  by National Natural Science Foundation of China. }}

\author{Li-Ping Huang\footnote{E-mail address: lipingmath@163.com (L. Huang)}\\
{\footnotesize  \em  School of Math. and Statis.,
Changsha University of Science and
Technology, Changsha, 410004, China}}
 \date{}
 \maketitle

\begin{abstract}
We investigate the generalized bilinear forms graph $\Gamma_d$ over a residue class ring $\mathbb{Z}_{p^s}$. We show that  $\Gamma_d$ is a connected vertex transitive graph,
and  completely determine its independence number, clique number, chromatic number and  maximum cliques. We also prove
that cores of both $\Gamma_d$ and its complement are maximum cliques. The graph $\Gamma_d$ is useful for error-correcting codes.
We show that every largest independent set of  $\Gamma_d$ is both an MRD code over $\mathbb{Z}_{p^s}$ and a usual MDS code.
Moreover, there is a largest independent set of $\Gamma_d$ to be a linear code over $\mathbb{Z}_{p^s}$.

\vspace{2mm}


\vspace{2mm}

\noindent{\bf 2010 AMS Classification}:   05C25, 15B33, 05C30, 94B60, 94B65

\end{abstract}

\section{Introduction}
 \setcounter{equation}{0}

\ \ \ \ \
Throughout, let $R$ be a commutative local ring and  $R^*$  the set of all units of $R$.  For a subset $S$ of $R$, let $S^{m\times n}$ be the set of all $m\times n$ matrices over $S$,
and let $S^n=S^{1\times n}$.
Let  $GL_n(R)$ be the set of $n\times n$ invertible matrices over $R$.
 Let $^tA$ denote the transpose matrix of a matrix $A$.
Denote by $I_r$ ($I$ for short)  the $r\times r$ identity matrix, and
${\rm diag}(A_1,\ldots, A_k)$ a block diagonal matrix where $A_i$ is an $m_i\times n_i$ matrix.  The cardinality of a set $X$ is denote by $|X|$.

For $0\neq A\in R^{m\times n}$, by Cohn's definition \cite{cohn2}, the {\em inner rank} of $A$,  denoted by $\rho(A)$, is the least  positive integer $r$ such that
\begin{equation}\label{innerrank1}
\mbox{$A=BC$ \  where  $B\in R^{m\times r}$ and $C\in R^{r\times n}$.}
\end{equation}
 Let $\rho(0)=0$.
For $A\in R^{m\times n}$, it is clear that $\rho(A)\leq \min\left\{m,n\right\}$ and $\rho(A)=0$ if and only if $A=0$.
When $R$ is a field, we have $\rho(A)={\rm rank}(A)$, where ${\rm rank}(A)$ is the usual rank of matrix over a field.
For matrices over $R$, we have (cf. \cite[Section 5.4]{cohn2,Freeideal}):
$\rho(A)=\rho(PAQ)$  where $P$ and $Q$ are invertible matrices over $R$; $\rho(A+B)\leq \rho(A)+\rho(B)$ and
$\rho(AC)\leq \min\left\{\rho(A), \rho(C) \right\}$.

For $A,B\in  R^{m\times n}$, the {\em rank distance} between $A$ and $B$ is defined by
\begin{equation}\label{rankdistance}
{\rm d_R}(A,B)=\rho(A-B).
\end{equation}
We have that $ {\rm d_R}(A, B)=0 \Leftrightarrow A=B$, ${\rm d_R}(A, B)={\rm d_R}(B, A)$ and ${\rm d_R}(A, B)\leq {\rm d_R}(A, C)+{\rm d_R}(C, B)$,
for all matrices of  appropriate sizes $A,B,C$ over $R$.

Let  $\mathbb{Z}_{p^s}$ denote the {\em residue class ring} of integers modulo $p^s$, where $p$ is a prime and  $s$ is a positive integer.
 The  $\mathbb{Z}_{p^s}$ is a {\em Galois ring}, a commutative {\em local ring},  a finite {\em principal ideal ring} (cf. \cite{mcdonald1, GaloisRing}).
  The principal ideal $(p)$ is the unique maximal ideal  of $\mathbb{Z}_{p^s}$,
 and denoted by $J_{p^s}$. The $J_{p^s}$ is also the {\em Jacbson radical} of $\mathbb{Z}_{p^s}$.
When $s=1$,  $\mathbb{Z}_{p}$ is a finite field with $p$ elements.
We have (cf. \cite{mcdonald1,GaloisRing}) that
\begin{equation}\label{d13fhmm8684}
\left|\mathbb{Z}_{p^s}\right|=p^s, \ \  \left|\mathbb{Z}_{p^s}^*\right|=(p-1)p^{s-1}, \ \ \left|J_{p^s}\right|=p^{s-1}.
\end{equation}

Let $\mathbb{F}_q$ be the finite field with $q$ elements (where $q$ is a power of a prime).
 All  graphs are {\em simple} \cite{Godsil} and finite in this paper.  Let $V(G)$ denote the vertex set of a graph $G$.
For $x,y\in V(G)$, we write $x \sim y$ if vertices $x$ and $y$  are adjacent.
Denote by ${\rm Aut}(G)$  the automorphism group  of a graph $G$.

The {\em generalized bilinear forms graph} over $\mathbb{F}_q$, denoted by  $\Gamma_d(\mathbb{F}_{q}^{m\times n})$,
has the vertex set $\mathbb{F}_{q}^{m\times n}$ where $m, n\geq 2$, and two distinct
vertices $A$ and $B$ are adjacent if ${\rm rank}(A-B)<d$ where $d$ is fixed with $2\leq d\leq {\rm min}\{m,n\}$.
When $d=2$, $\Gamma_2(\mathbb{F}_{q}^{m\times n})$ is  the usual {\em bilinear forms graph} over $\mathbb{F}_q$.
The  bilinear forms graph plays an important role in combinatorics and coding theory, and
it has been extensively studied (cf. \cite{Brouwera2, Delsarte,  Gadouleau,Huang-2014, Tanaka, Y.X.Wang}).

Recently, the bilinear forms graph over  $\mathbb{Z}_{p^s}$ is studied by \cite{Huang-2017}. However, the generalized bilinear forms graph over  $\mathbb{Z}_{p^s}$
 remains to be further studied.  As a natural extension of the generalized bilinear forms graph over  $\mathbb{F}_q$,
we define the generalized bilinear forms graph over  $\mathbb{Z}_{p^s}$ as follows.
The {\em generalized bilinear forms graph} ({\em bilinear forms graph} for short) over $\mathbb{Z}_{p^s}$, denoted by  $\Gamma_d(\mathbb{Z}_{p^s}^{m\times n})$
($\Gamma_d$ for short), has the vertex set $\mathbb{Z}_{p^s}^{m\times n}$ where $m,n\geq 2$, and two distinct
vertices $A$ and $B$ are adjacent if $\rho(A-B)<d$, where $d$ is fixed with $2\leq d\leq {\rm min}\{m,n\}$.

MRD codes  and codes over $\mathbb{Z}_{p^s}$  are active research topics in the coding theory (cf. \cite{Delsarte, Gabidulin, Roth,Cossidente, Horlemann, Neri})
and \cite{Dougherty2,Oued,Lee}). The generalized bilinear forms graph $\Gamma_d$ has good  application to the error-correcting codes over $\mathbb{Z}_{p^s}$. In fact,
we will show that every largest independent set of  $\Gamma_d$ is both an  MRD code over $\mathbb{Z}_{p^s}$ and a usual MDS code \cite{Roth2}.
Moreover, there is a largest independent set of $\Gamma_d$ such that it is a linear code over $\mathbb{Z}_{p^s}$.

The  paper is organized as follows.  In Section 2, we recall some properties of matrices over $\mathbb{Z}_{p^s}$.
In Section 3,  we show that  $\Gamma_d$ is a connected vertex transitive graph, and determine  the independence number,  the clique number and the chromatic number of $\Gamma_d$. We also show that
every largest independent set of  $\Gamma_d$ is both an $(m\times n, d)$ MRD code over $\mathbb{Z}_{p^s}$ and a usual MDS code, and
there is a largest independent set of $\Gamma_d$ to be a linear code over $\mathbb{Z}_{p^s}$.
In Section 4, We will determine the algebraic structures of maximum cliques of $\Gamma_d$, and
  show that cores of  both $\Gamma_d$ and its complement are maximum cliques.

\section{Matrices  over $\mathbb{Z}_{p^s}$}
 \setcounter{equation}{0}

\ \ \ \ \  In this section,  we recall some basic properties of matrices over $\mathbb{Z}_{p^s}$.

\begin{lem}\label{uniqueuk1}{\rm (see \cite[Proposition 6.2.2]{Bini} or \cite[p.328]{mcdonald1})} \ Every non-zero element $x$ in $\mathbb{Z}_{p^s}$ can be written as
$x=up^t$ where $u$ is a unit and  $0\leq t\leq s-1$. Moreover, the integer $t$ is unique and $u$ is unique modulo the ideal $(p^{s-t})$ of $\mathbb{Z}_{p^s}$.
\end{lem}

Let $T_p=\{0,1,\ldots, p-1\}\subseteq\mathbb{Z}_{p^s}$.
For two distinct elements $a,b\in T_p$,  we always have $a-b\in \mathbb{Z}_{p^s}^*$.
Without loss of generality, we may assume that $T_p=\mathbb{Z}_p$ in our discussion.

\begin{lem}\label{uniqueuk2}{\rm (cf. \cite[p.328]{mcdonald1})} \  Every non-zero element $x$ in $\mathbb{Z}_{p^s}$ can be written uniquely as
$$x=t_0+t_1p+\cdots +t_{s-1}p^{s-1},$$
where $t_i\in T_p$,  $i=0,1,\ldots, s-1$.
\end{lem}

By Lemma \ref{uniqueuk2},
every matrix $X\in\mathbb{Z}_{p^s}^{m\times n}$ can be written uniquely as
\begin{equation}\label{bcweu7754332t}
X=X_0+X_1p+\cdots +X_{s-1}p^{s-1},
\end{equation}
where $X_i\in T_p^{m\times n}$, $i=0,\ldots, s-1$.

Note that every matrix in $T_p^{m\times n}$ can be seen as a matrix in $\mathbb{Z}_p^{m\times n}$.
We define the {\em natural surjection}
\begin{equation}\label{naturalsurjection}
\pi: \mathbb{Z}_{p^s}^{m\times n}\rightarrow \mathbb{Z}_{p}^{m\times n}
\end{equation}
by $\pi(X)=X_0$ for all $X\in\mathbb{Z}_{p^s}^{m\times n}$ of the form (\ref{bcweu7754332t}). Clearly, $\pi(A)=A$ if $A\in\mathbb{Z}_p^{m\times n}$.
For  $X, Y\in\mathbb{Z}_{p^s}^{m\times n}$ and $Q\in\mathbb{Z}_{p^s}^{n\times k}$, We have
\begin{equation}\label{ytr64xceqw13ityi}
 \pi(X+Y)=\pi(X)+\pi(Y),
 \end{equation}
\begin{equation}\label{pipipipi003}
 \pi(XQ)=\pi(X)\pi(Q),
\end{equation}
\begin{equation}\label{pipipipi004}
 \pi(^tX)= \,^t(\pi(X)).
\end{equation}

By Lemma \ref{uniqueuk1}, it is easy to prove the following result.

\begin{lem}\label{canonical}
{\rm (cf. \cite[Chap. II]{Newman} or \cite[p.327]{mcdonald1})} \
Let $R=\mathbb{Z}_{p^s}$ where $s\geq 2$, and let $A\in R^{m\times n}$ be a non-zero matrix.
Then there are $P\in GL_m(R)$ and $Q\in GL_n(R)$ such that
\begin{equation}\label{form1}
A=P {\rm diag}\left(I_r, p^{k_1}, \ldots, p^{k_t}, 0  \right)Q,
\end{equation}
where  $1\leq k_1\leq \cdots\leq k_t\leq s-1$. Moreover, the parameters $(r, t, k_1,\ldots, k_t)$ are uniquely determined by $A$.
In (\ref{form1}), $I_r$ or ${\rm diag}\left(p^{k_1},\ldots,p^{k_t}\right)$ may be absent.
\end{lem}

Let  $A\in R^{m\times n}$, and let  $I_k(A)$ be the ideal in $R$ generated by all $k\times k$ minors of $A$, $k=1,\ldots, {\rm min}\{m,n\}$.
Let ${\rm Ann}_R(I_k(A))=\left\{x\in R: xI_k(A)=0\right\}$ denote the annihilator of $I_k(A)$. The {\em McCoy rank} of $A$, denoted by ${\rm rk}(A)$, is the following integer:
$${\rm rk}(A)={\rm max}\left\{k: {\rm Ann}_R(I_k(A))=(0)\right\}.$$
 We have that ${\rm rk}(A)={\rm rk}(^tA)$;
${\rm rk}(A)={\rm rk}(PAQ)$ where $P$ and $Q$ are invertible matrices of the appropriate sizes; and ${\rm rk}(A)=0$ if and only if ${\rm Ann}_R(I_1(A))\neq (0)$ (cf. \cite{brown}).

\begin{lem}\label{innerrank02} {\rm(see \cite[Lemmas 2.4 and 2.7]{Huang-2017})} \
Let $0\neq A\in \mathbb{Z}_{p^s}^{m\times n}$ ($s\geq 2$) be of the form (\ref{form1}). Then $r+t$ is  the inner rank of $A$,
and $r$ is  the McCoy rank of $A$.
\end{lem}

Let $A\in \mathbb{Z}_{p^s}^{m\times n}$. By Lemma \ref{innerrank02}, ${\rm rk}(A)\leq \rho(A)$, and
${\rm rk}(A)=0$ if and only if $A\in  J_{p^s}^{m\times n}$.
 By Lemmas \ref{canonical} and \ref{innerrank02}, it is easy to see that (cf. \cite{Huang-2017})
\begin{equation}\label{rank-0008}
\mbox{$\rho\left(
            \begin{array}{cc}
              A & 0 \\
              0 & B\\
            \end{array}
          \right)=\rho(A)+\rho(B)$, \ \ if $A,B$ are matrices over $\mathbb{Z}_{p^s}$.}
\end{equation}

For $A\in R^{m\times n}$ and  $B\in R^{n\times m}$, if $AB=I_m$, we call that $A$ has a {\em right inverse} and $B$ is a right inverse of $A$.
Similarly, if $AB=I_m$, than  $B$ has a {\em left inverse} and $A$ is a left inverse of $B$.  Note that if $a\in\mathbb{Z}_{p^s}^*$ and $b\in J_{p^s}$,
then $a\pm b\in\mathbb{Z}_{p^s}^*$. By Lemma \ref{canonical}, we have the following lemmas.

\begin{lem}\label{Mc-rank-1}{\rm (cf. \cite{Huang-2017})} \
Let $A\in \mathbb{Z}_{p^s}^{m\times n}$ where $s\geq 2$ and $n\geq m$. Then ${\rm rk}(A)={\rm rk}(A\pm B)$ for all $B\in J_{p^s}^{m\times n}$.
Moreover, $A$ has a right inverse  if and only if ${\rm rk}(A)=m$.
\end{lem}

\begin{lem}\label{ds424rrwrww}{\rm (see. \cite[Lemma 4.2]{Huang-2017})} \  If $A\in GL_n(\mathbb{Z}_{p^{s-1}})$ where $s\geq 2$, then $A\in GL_n(\mathbb{Z}_{p^{s}})$.
 \end{lem}

\begin{lem}\label{64gdgd7wrww}{\rm (see \cite[Lemma 4.3]{Huang-2017})} \  If  $A\in\mathbb{Z}_{p^{s-1}}^{m\times n}$ where $s\geq 2$, then both $A$ and $Ap$
can be viewed as matrices in $\mathbb{Z}_{p^{s}}^{m\times n}$ with the same inner rank.
\end{lem}

 \section{Independence number of $\Gamma_d$ and MDR codes}
 \setcounter{equation}{0}

\subsection{Independence number and chromatic number of $\Gamma_d$}

\ \ \ \ \
Recall that an {\em independent set} of  a graph $G$ is a subset of vertices such that no two vertices are adjacent. A {\em largest independent set}
of $G$ is an independent set of maximum cardinality. The {\em independence number} of $G$, denoted by $\alpha(G)$,  is  the number of vertices in a largest independent set of $G$.

 An {\em $r$-colouring} of a graph $G$ is a homomorphism from $G$ to  the complete graph $K_r$.
The {\em chromatic number} of $G$, denoted by $\chi(G)$, is  the least value $k$ for which $G$ can be  $k$-coloured.

A {\em clique} of a graph $G$ is a complete subgraph of $G$. A clique $C$ is  maximal if there is no clique of $G$ which properly contains $C$ as a subset.
A {\em maximum clique} of $G$ is a clique of $G$ which has maximum cardinality.
The  clique number of $G$, denoted by $\omega(G)$, is the number of vertices in a  maximum clique.
For convenience, we regard that a maximal clique and its vertex set are the same.

\begin{thm}\label{dfs424rew}  \ Every generalized bilinear forms graph $\Gamma_d(\mathbb{Z}_{p^s}^{m\times n})$ is a connected vertex transitive graph.
\end{thm}
\proof
Let $R=\mathbb{Z}_{p^{s}}$. For any vertex $A$ of $\Gamma_d$,
since the  map $X\mapsto X-A$ is an automorphism of  $\Gamma_d$,  $\Gamma_d$ is vertex-transitive.

 Let $A, B\in R^{m\times n}$ with $\rho(A-B)=r>0$. By Lemma \ref{canonical},
 there are $P\in GL_m(R)$ and $Q\in GL_n(R)$ such that
$B-A=P{\rm diag}\left(p^{k_1}, \ldots, p^{k_r}, 0,\ldots, 0 \right)Q$, where  $0\leq k_1\leq \cdots\leq k_r\leq {\rm max}\{s-1, 1\}$. If $r\leq d-1$, then $A\sim B$.
Now we assume that $r>d-1$.
Put $A_0=A$, $A_i=A+P{\rm diag}\left(p^{k_1}, \ldots, p^{k_i}, 0,\ldots, 0 \right)Q$, $i=1,\ldots,r$. Then $A_r=B$ and
$A_i\sim A_{i+1}$, $i=0,1,\ldots, r-1$. It follows that $\Gamma_d$ is  connected.
$\qed$

For a graph $G$,  we have (see \cite[Theorem 6.10, Corollary 6.2]{Chartrand})
\begin{equation}\label{chromatic-1}
 \chi( G ) \geq  {\rm max} \left\{\omega( G ), \ |V( G )|/\alpha( G ) \right\}.
\end{equation}

From \cite[Lemma 2.7.2]{Roberson} we have
\begin{equation}\label{d34223gh}
\chi(G)\geq \frac{|V(G)|}{\alpha(G)}\geq \omega(G), \ \mbox{if $G$ is  vertex-transitive.}
\end{equation}

 Now, we recall the coding theory on a finite field $\mathbb{F}_{q}$ (cf. \cite{Cossidente, Gadouleau, Horlemann, Neri}). Without loss of generality
 we assume that  $n\geq m\geq d >1$ are integers.
 An  $(m\times n, d)$ {\em rank distance code} over $\mathbb{F}_{q}$ is a subset  $\mathcal{C}$ of $\mathbb{F}_{q}^{m\times n}$
with  ${\rm d_R}(A,B)\geq d$ for distinct $A,B\in \mathcal{C}$.
For an $(m\times n, d)$ rank distance code $\mathcal{C}$ over $\mathbb{F}_{q}$, we have $|\mathcal{C}|\leq q^{n(m-d+1)}$, and
the bound $q^{n(m-d+1)}$ is called the {\em Singleton bound} for $\mathcal{C}$.
If a rank distance code $\mathcal{C}$ over $\mathbb{F}_{q}$ satisfies $|\mathcal{C}|=q^{n(m-d+1)}$, then  $\mathcal{C}$ is called
an $(m\times n, d)$ {\em maximum rank distance code} (MRD code in short) over $\mathbb{F}_{q}$.
A {\em linear code} of length $n'$ over $\mathbb{F}_{q}$ is a subspace of $\mathbb{F}_{q}^{n'}$.

  MRD codes (over $\mathbb{F}_{q}$) can be used to correct errors and erasures in network. In 1978,
Delsarte \cite{Delsarte} (and independently  Gabidulin in 1985 \cite{Gabidulin},  Roth in 1991 \cite{Roth}) proved the following important result:

\begin{lem}\label{MRDcodes01} \
There are linear $(m\times n, d)$ MRD codes over $\mathbb{F}_{q}$ for all choices of $m,n,d$.
\end{lem}

These linear MRD codes in Lemma \ref{MRDcodes01} are also called {\em Gabidulin codes}. Until a few years ago, the only known MRD codes were Gabidulin codes.
Recently, the research of MRD codes is active, and we know about other some MRD codes over $\mathbb{F}_{q}$ (cf. \cite{Cossidente, Horlemann, Neri}).

 Clearly, every independent set of $\Gamma_d(\mathbb{F}_{q}^{m\times n})$ is an $(m\times n, d)$ rank distance code over $\mathbb{F}_{q}$ and vice versa.
 In other words, an $(m\times n, d)$ MRD code over $\mathbb{F}_{q}$ is a largest independent set of $\Gamma_d(\mathbb{F}_{q}^{m\times n})$.
 Thus, Lemma \ref{MRDcodes01} implies that
 \begin{equation}\label{gene4fsrfsfs}
\mbox{$\alpha\left( \Gamma_d(\mathbb{F}_{q}^{m\times n})\right)=q^{n(m-d+1)}$ (where $n\geq m$).}
\end{equation}
The formula (\ref{gene4fsrfsfs}) was also showed by  \cite{Gadouleau}.

By (\ref{d34223gh}) and (\ref{gene4fsrfsfs}), we get $\omega\left( \Gamma_d(\mathbb{F}_{q}^{m\times n})\right)\leq q^{n(d-1)}$. On the other hand,
 $\small\mathcal{M}:=\left(
  \begin{array}{c}
   \mathbb{F}_{q}^{(d-1)\times n} \\
    0 \\
     \end{array}
   \right)$ is a clique and $|\mathcal{M}|=q^{n(d-1)}$. Thus we obtain
\begin{equation}\label{geneytr675vxcwq128}
\mbox{$\omega\left( \Gamma_d(\mathbb{F}_{q}^{m\times n})\right)=q^{n(d-1)}$ (where $n\geq m$).}
\end{equation}

\begin{thm}\label{bilinearformsyr56bch8} \ When $n\geq m$, the  independence number  and the clique number of $\Gamma_d(\mathbb{Z}_{p^s}^{m\times n})$ are
\begin{equation}\label{bilinearform5435gdfgd865}
\alpha\left(\Gamma_d(\mathbb{Z}_{p^s}^{m\times n})\right)= p^{sn(m-d+1)}
\end{equation}
and
\begin{equation}\label{biline13nvhiu6}
\omega\left(\Gamma_d(\mathbb{Z}_{p^s}^{m\times n})\right)=p^{sn(d-1)}.
\end{equation}
\end{thm}
\proof
We prove (\ref{bilinearform5435gdfgd865}) and (\ref{biline13nvhiu6}) by induction on $s$. When $s=1$, (\ref{bilinearform5435gdfgd865}) and (\ref{biline13nvhiu6})
hold by (\ref{gene4fsrfsfs}) and (\ref{geneytr675vxcwq128}). Suppose that $s\geq 2$ and
$$
\alpha\left(\Gamma_d(\mathbb{Z}_{p^{s-1}}^{m\times n})\right)= p^{(s-1)n(m-d+1)}.
$$
Write $\alpha=p^{n(m-d+1)}$. Let $\mathcal{A}_1=\left\{S_1, S_2, \ldots, S_{\alpha}\right\}$ be a largest independent set
of $\Gamma_d(\mathbb{Z}_{p}^{m\times n})$, where $S_i\in\mathbb{Z}_{p}^{m\times n}$, $i=1, \ldots, \alpha$.
Since each $S_i$ can be seen as an element in $\mathbb{Z}_{p^s}^{m\times n}$,
 $\mathcal{A}_1$ is also an independent set of $\Gamma_d(\mathbb{Z}_{p^s}^{m\times n})$.

Put $\beta= p^{(s-1)n(m-d+1)}$.
Let $\mathcal{A}_{s-1}=\{C_1,C_2,\ldots,C_\beta\}$ be a largest independent set of  $\Gamma_d(\mathbb{Z}_{p^{s-1}}^{m\times n})$.
Note that every matrix over $\mathbb{Z}_{p^{s-1}}$ can be seen as a matrix over $\mathbb{Z}_{p^{s}}$.
For any two distinct matrices $C_i, C_j\in\mathcal{A}_{s-1}$, we have  $\rho(C_i-C_j)\geq d$.  From Lemma \ref{64gdgd7wrww} we get  $\rho(C_ip-C_jp)\geq d$.
Thus $\mathcal{A}_{s-1}p:=\{Xp: X\in\mathcal{A}_{s-1}\}$ is  an independent set of  $\Gamma_d(\mathbb{Z}_{p^{s}}^{m\times n})$  with $|\mathcal{A}_{s-1}p|=\beta$.
Let
\begin{equation}\label{64gsktgh7ghfk75h}
 \mathcal{L}_k=S_k+\mathcal{A}_{s-1}p=\left\{S_k+C_1p, S_k+C_2p, \ldots, S_k+C_\beta p \right\}, \ \ k=1,\ldots, \alpha.
 \end{equation}
Then
$\mathcal{L}_k$ is an independent set of  $\Gamma_d(\mathbb{Z}_{p^s}^{m\times n})$ and $|\mathcal{L}_k|=\beta$, $k=1,\ldots, \alpha$.
Let
\begin{equation}\label{dseeqt7ijghbd4y}
\mathcal{A}_s=\mathcal{L}_1\cup \mathcal{L}_2\cup\cdots \cup \mathcal{L}_{\alpha}:=\mathcal{A}_1+\mathcal{A}_{s-1}p.
\end{equation}
Since ${\rm rk}(S_i-S_j)\geq d$ for $i\neq j$,
 Lemma \ref{Mc-rank-1} implies that  $\mathcal{A}_s$ is an independent set of  $\Gamma_d(\mathbb{Z}_{p^s}^{m\times n})$ with $|\mathcal{A}_s|=\alpha\beta=p^{sn(m-d+1)}$.
Therefore,
\begin{equation}\label{gf423jgkq9679bv}
\alpha\left(\Gamma_d(\mathbb{Z}_{p^s}^{m\times n})\right)\geq p^{sn(m-d+1)}.
\end{equation}

By (\ref{d34223gh}) and (\ref{gf423jgkq9679bv}), we have
$$\omega\left(\Gamma_d(\mathbb{Z}_{p^{s}}^{m\times n})\right)\leq \frac{\left|V\left(\Gamma_d(\mathbb{Z}_{p^{s}}^{m\times n}\right)\right|}
{\alpha\left(\Gamma_d(\mathbb{Z}_{p^{s}}^{m\times n})\right)}
\leq \frac{p^{smn}}{p^{sn(m-d+1)}}= p^{sn(d-1)}.$$
On the other hand,  it is easy to see that  $\small\mathcal{M}:=\left(
                                      \begin{array}{c}
                                        \mathbb{Z}_{p^s}^{(d-1)\times n} \\
                                        0 \\
                                      \end{array}
                                    \right)$ is a clique of  $\Gamma_d(\mathbb{Z}_{p^s}^{m\times n})$ and $|\mathcal{M}|=p^{sn(d-1)}$. Thus we obtain
\begin{equation}\label{fd235fghhjgalju0}
\omega\left(\Gamma_d(\mathbb{Z}_{p^s}^{m\times n})\right)= p^{sn(d-1)}.
\end{equation}

Using (\ref{d34223gh}) and (\ref{fd235fghhjgalju0}), we get
$$\alpha\left(\Gamma_d(\mathbb{Z}_{p^{s}}^{m\times n})\right)\leq \frac{\left|V\left(\Gamma_d(\mathbb{Z}_{p^{s}}^{m\times n}\right)\right|}
{\omega\left(\Gamma_d(\mathbb{Z}_{p^{s}}^{m\times n})\right)}=\frac{p^{smn}}{p^{sn(d-1)}}= p^{sn(m-d+1)}.$$
It follows from (\ref{gf423jgkq9679bv}) that (\ref{bilinearform5435gdfgd865}) holds.
 $\qed$

Let $G$ be a finite group and let $C$ be a subset of $G$ that is closed under taking inverses and does not contain the identity. The  {\em Cayley graph} $X(G,C)$
is the graph with vertex set $G$ and two vertices $h$ and $g$ are adjacent if  $hg^{-1}\in C$. A Cayley graph $X(G,C)$ is {\em normal} if
$gCg^{-1}=C$ for all $g\in G$.
It is easy to see that $\Gamma_d(\mathbb{Z}_{p^s}^{m\times n})$ is a normal Cayley graph on  the matrix additive group $G$ of $\mathbb{Z}_{p^s}^{m\times n}$
and the inverse closed subset $C$  is the set of all matrices of inner rank $<d$.

\begin{lem}\label{hf5346hfhfhfss}{\rm (Godsil \cite[Corollary 6.1.3]{Godsil3}, cf. \cite[Therem 3.3.1]{Roberson})} \
Let $G$ be a normal Cayley graph. If $\alpha(G)\omega(G)=|V(G)|$, then $\chi(G)=\omega(G)$.
\end{lem}

By  Theorem \ref{bilinearformsyr56bch8} and Lemma \ref{hf5346hfhfhfss}, the chromatic number of $\Gamma_d(\mathbb{Z}_{p^s}^{m\times n})$ is

\begin{equation}\label{gf353tgdcvadkhl54}
\mbox{$\chi\left(\Gamma_d(\mathbb{Z}_{p^s}^{m\times n})\right)=\omega\left(\Gamma_d(\mathbb{Z}_{p^s}^{m\times n})\right)=p^{sn(d-1)}$ (where $n\geq m$)} .
\end{equation}

\subsection{ MDR codes over  $\mathbb{Z}_{p^s}$}

\ \ \ \ \
We recall the usual definition of MDS code in coding theory (cf. \cite{Roth2}).
An $(n,M)$ {\em code} over a finite alphabet $F$ is a nonempty subset $\mathcal{C}$ of size $M$ of $F^n$.
For  codewords $\alpha=(a_1,\ldots, a_n)$, $\beta=(b_1,\ldots, b_n)$, the {\em Hamming distance} between $\alpha$ and $\beta$ is
$d(\alpha,\beta) = \left|\{i: a_i \neq b_i, 1\leq i\leq n \}\right|$.
An $(n,M)$ code with minimum distance $d$ is called an {\em $(n,M, d)$ code}. For any $(n,M, d)$ code over an alphabet of size $q$, we have the {\em Singleton bound}
$$d\leq n-({\rm log}_qM)+1.$$
 An $(n,M, d)$ code over an alphabet of size $q$ is called a {\em maximum distance separable code}
(MDS code in short) if it attains the Singleton bound, i.e. $d= n-({\rm log}_qM)+1$.

A {\em code} $\mathcal{C}$ (over $\mathbb{Z}_{p^s}$) of length $m$ is a subset of $\mathbb{Z}_{p^s}^{m\times 1}$ (or $\mathbb{Z}_{p^s}^{m}$). 
If the code is a submodule (i.e. vector space) over $\mathbb{Z}_{p^s}$
we say that it is a {\em linear code} over $\mathbb{Z}_{p^s}$.
Suppose that  $n\geq m\geq d >1$ are integers.
 As a natural extension of rank distance code over $\mathbb{F}_{q}$, we define the rank distance code over $\mathbb{Z}_{p^s}$ as follows.
 An  $(m\times n, d)$ {\em rank distance code} over $\mathbb{Z}_{p^s}$ is a subset  $\mathcal{C}$ of $\mathbb{Z}_{p^s}^{m\times n}$
with  ${\rm d_R}(A,B)\geq d$ for distinct $A,B\in \mathcal{C}$.

 Every independent set of $\Gamma_d(\mathbb{Z}_{p^s}^{m\times n})$ is
an $(m\times n, d)$ rank distance code over $\mathbb{Z}_{p^s}$ and vice versa.
Thus, for an $(m\times n, d)$ rank distance code $\mathcal{C}$ over $\mathbb{Z}_{p^s}$, we have $|\mathcal{C}|\leq p^{sn(m-d+1)}$ by (\ref{bilinearform5435gdfgd865}).
If a rank distance code $\mathcal{C}$ over $\mathbb{Z}_{p^s}$ satisfies $|\mathcal{C}|=q^{n(m-d+1)}$, then  $\mathcal{C}$
is called an $(m\times n, d)$ {\em maximum rank distance code} (MRD code in short) over $\mathbb{Z}_{p^s}$.
In other words, an $(m\times n, d)$ MRD code over $\mathbb{Z}_{p^s}$ is a largest independent set of  $\Gamma_d(\mathbb{Z}_{p^{s}}^{m\times n})$.

\begin{thm}\label{MRDcodes02} \ If $\mathcal{S}$ is a largest independent set of  $\Gamma_d(\mathbb{Z}_{p^{s}}^{m\times n})$, then $\mathcal{S}$
 is both an $(m\times n, d)$ MRD code over $\mathbb{Z}_{p^s}$ and an $(m,M, d)$ MDS code  with $M=p^{sn(m-d+1)}$.
 Moreover,  there is a largest independent set  $\mathcal{C}$ of  $\Gamma_d(\mathbb{Z}_{p^{s}}^{m\times n})$ such that $\mathcal{C}$ is
a linear code over $\mathbb{Z}_{p^s}$.
\end{thm}
\proof
{\bf Step 1.} \ Without loss of generality we assume that  $n\geq m\geq d >1$ are integers.
By  (\ref{bilinearform5435gdfgd865}), it is clear that every largest independent set of $\Gamma_d(\mathbb{Z}_{p^{s}}^{m\times n})$ is
 an $(m\times n, d)$ MRD code over $\mathbb{Z}_{p^s}$. We prove that every largest independent set of $\Gamma_d(\mathbb{Z}_{p^{s}}^{m\times n})$ is
 a usual MDS code as follows.

By Lemma \ref{uniqueuk2} and (\ref{bcweu7754332t}),  the row vector space $\mathbb{Z}_{p^{s}}^{n}$  is isomorphic to $\mathbb{Z}_{p^{sn}}$
(as an $n$-dimensional vector space over $\mathbb{Z}_{p^{s}}$) (cf. \cite[Chapter 14]{GaloisRing}).
Thus, it is easy to see that $\mathbb{Z}_{p^{s}}^{m\times n}$  (as a vector space over $\mathbb{Z}_{p^{s}}$) is isomorphic
to the column vector space $\mathbb{Z}_{p^{sn}}^{m\times 1}$ (as a vector space over $\mathbb{Z}_{p^{s}}$).
Let $\mathcal{S}$ be a largest independent set of  $\Gamma_d(\mathbb{Z}_{p^{s}}^{m\times n})$. From (\ref{bilinearform5435gdfgd865}) we have 
$|\mathcal{S}|=p^{sn(m-d+1)}$.
 For two distinct vertices  $A,B\in\mathcal{S}$, since $\rho(A-B)\geq d$, it is clear that $A-B$ has at least $d$ non-zero rows.
 Hence the Hamming distance between $A$ and $B$ (as vectors in
$\mathbb{Z}_{p^{sn}}^{m\times 1}$) is at least $d$. Thus, $\mathcal{S}$ can be seen as an $(m, M, d)$ code
over an alphabet of size $p^{sn}$, where $M=|\mathcal{S}|=p^{sn(m-d+1)}$. Since $m-({\rm log}_{p^{sn}}M)+1=d$, $\mathcal{S}$ is a usual MDS code.

 {\bf Step 2.} \ We assert that there is a largest independent set  $\mathcal{C}$ of  $\Gamma_d(\mathbb{Z}_{p^{s}}^{m\times n})$ such that it is
a linear codes over $\mathbb{Z}_{p^{s}}$. By Lemma \ref{MRDcodes01}, we  may assume with no loss of generality that $s\geq 2$. Using Lemma \ref{MRDcodes01},
there is a linear $(m\times n, d)$ MRD code $\mathcal{A}_1$ over $\mathbb{Z}_{p}$. Clearly, $\mathcal{A}_1$ is a vector space over $\mathbb{Z}_{p}$ and
a largest independent set of  $\Gamma_d(\mathbb{Z}_{p}^{m\times n})$.

When $s=2$, recalling  (\ref{dseeqt7ijghbd4y}),
$\mathcal{A}_2:=\mathcal{A}_1+\mathcal{A}_1p$ is a largest independent set of  $\Gamma_d(\mathbb{Z}_{p^2}^{m\times n})$.
Since $\mathcal{A}_1$ contains $0$, we have $\mathcal{A}_1\subset\mathcal{A}_{2}$. By Lemma \ref{uniqueuk2} and (\ref{bcweu7754332t}),
 one can prove that $\mathcal{A}_2$ is a vector space over $\mathbb{Z}_{p^2}$.
For $s-1$ ($s\geq 3$), assume that there is a
largest independent set  $\mathcal{A}_{s-1}$ of $\Gamma_d(\mathbb{Z}_{p^{s-1}}^{m\times n})$, such that $\mathcal{A}_1\subset\mathcal{A}_{s-1}$
and $\mathcal{A}_{s-1}$ is a vector space over $\mathbb{Z}_{p^{s-1}}$.
Then by (\ref{dseeqt7ijghbd4y}), $\mathcal{A}_s:=\mathcal{A}_1+\mathcal{A}_{s-1}p$ is a largest independent set of  $\Gamma_d(\mathbb{Z}_{p^s}^{m\times n})$.
 Since $\mathcal{A}_{s-1}$ contains $0$,  $\mathcal{A}_1\subset\mathcal{A}_{s}$. Applying Lemma \ref{uniqueuk2} and (\ref{bcweu7754332t}), it is easy to prove that
 $\mathcal{A}_{s}$ is a vector space over $\mathbb{Z}_{p^s}$.
 By the induction on $s$, for any  $s\geq 2$, there exists a largest independent set  $\mathcal{A}_s:=\mathcal{A}_1+\mathcal{A}_{s-1}p$ of  $\Gamma_d(\mathbb{Z}_{p^s}^{m\times n})$,
 such that $\mathcal{A}_s$ is a vector space over $\mathbb{Z}_{p^s}$.
By Step 1,  $\mathcal{A}_s$ is a linear code over $\mathbb{Z}_{p^s}$.
$\qed$

By Theorem \ref{MRDcodes02}, there are linear $(m\times n, d)$ MRD codes over $\mathbb{Z}_{p^s}$ for all choices of $m,n,d$.
Applying  (\ref{dseeqt7ijghbd4y}), we can Construct many  MRD codes over $\mathbb{Z}_{p^s}$.

 \section{Maximum cliques  and Core of $\Gamma_d$}
 \setcounter{equation}{0}

\subsection{Maximum cliques of $\Gamma_d$}

\ \ \ \ \
Note that the algebraic structures of maximum cliques of $\Gamma_d$ have many applications. For example, a maximum clique of $\Gamma_d$
is a largest independent set of the complement of $\Gamma_d$.
We will determine the algebraic structures of maximum cliques of $\Gamma_d$.

Let $X=[\alpha_1,\ldots,\alpha_m]$ be an $m$-dimensional vector subspace of $\mathbb{F}_{q}^{n}$ $(n\geq m)$.
Then $X$ has a {\em matrix representation}
$\scriptsize \left(
     \begin{array}{c}
       \alpha_1 \\
       \vdots \\
       \alpha_m \\
     \end{array}
   \right)\in\mathbb{F}_{q}^{m\times n}$.
For simpleness, the matrix representation of a subspace $X$ is also denoted by $X$.
If $X$ is a matrix representation of a subspace $X$ of $\mathbb{F}_{q}^{n}$, then
$PX$ is also a matrix representation of $X$ where $P\in GL_m(\mathbb{F}_{q})$.
It follows that the matrix representation is not unique. However,
a subspace $X$ of $\mathbb{F}_{q}^{n}$ has a unique matrix representation which is the {\em row-reduced echelon form} $X=(I_m, B)Q$, where $Q$ is a permutation matrix.

A geometric description of the bilinear forms graph $\Gamma_2(\mathbb{F}_{q}^{m\times n})$ is the adjacency graph of the attenuated space.
Let  $n\geq m$ and let
$$W=(I_n,0)$$
 be a fixed $n$-dimension subspace of $\mathbb{F}_q^{m+n}$. Write
$$\mathcal{A}_{t}=\left\{U: U\subseteq\mathbb{F}_q^{m+n},  {\rm dim}(U)=t,  U\cap W=\{0\}\right\}, \ t=m-1,m.$$
The incidence structure $\left( \mathcal{A}_{m}, \mathcal{A}_{m-1}, \subseteq \right)$ is called an {\em attenuated space}. Its {\em adjacency graph} is  the graph
with $\mathcal{A}_{m}$ as its vertex set, and two vertices being adjacent if   their intersection is in $\mathcal{A}_{m-1}$.
For any $U\in\mathcal{A}_{m}$, by $U\cap W=\{0\}$ we have $U=(X_U,I_m)$, where $X_U\in\mathbb{F}_q^{m\times n}$ is uniquely determined by $U$.
Let
\begin{equation}\label{isomorphismattenuated}
\varphi(U)=X_U, \  U=(X_U,I_m)\in\mathcal{A}_{m}.
\end{equation}
 Then $\varphi$ is a graph isomorphism  from the adjacency graph of $\left( \mathcal{A}_{m-1}, \mathcal{A}_{m}, \subseteq \right)$ to
 $\Gamma_2(\mathbb{F}_{q}^{m\times n})$.
Moreover,   ${\rm dim}(U_1\cap U_2)=r$ if and only if
${\rm rank}(\varphi(U_1)-\varphi(U_2))=m-r$ for all $U_1,U_2\in \mathcal{A}_{m}$ (cf. \cite{T.Huang}, \cite[\S  9.5A]{Brouwera2}).

\begin{lem}\label{Rectangular-PID2-4a}{\em (see \cite[Theorem 3(2)]{Tanaka})} \
Suppose that $n\geq m$ and  $\Gamma$ is the adjacency graph of an attenuated space $\left( \mathcal{A}_{m}, \mathcal{A}_{m-1}, \subseteq \right)$.
Let $\mathscr{F}$ be a collection of elements of the $\mathcal{A}_{m}$ with the property that ${\rm dim}(\gamma \cap \delta)\geq t$
 for all $\gamma, \delta$ in $\mathscr{F}$. Then $|\mathscr{F}|\leq q^{(m-t)n}$, and equality holds if and only if either
  {\em (a)}  $\mathscr{F}$ consists of all elements of $\mathcal{A}_{m}$ which contain a fixed $t$-dimensional subspace $U$ with $U\cap W=0$, or
 {\em (b)} $m=n$ and $\mathscr{F}$ is the set of all elements of  $\mathcal{A}_{m}$ contained in a fixed  $(2m-t)$-dimensional subspace $U'$ with
 ${\rm dim}(U'\cap W)=m-t$.
\end{lem}

\begin{lem}\label{Rectangular-PID2-4} \  In $\Gamma_d(\mathbb{F}_{q}^{m\times n})$ (where $n\geq m$),  every maximum clique ${\cal M}$ containing $0$ is of  the
 form either
 \begin{equation}\label{maxiumclique35trgr}
 \mbox{${\cal M}=P\left(\begin{array}{c}
                              \mathbb{F}_{q}^{(d-1)\times n} \\
                              0 \\
                            \end{array}
                          \right):=\left\{P\left(
                            \begin{array}{c}
                              X\\
                              0 \\
                            \end{array}
                          \right): X\in\mathbb{F}_{q}^{(d-1)\times n} \right\}$,}
 \end{equation}
 where $P\in GL_m(\mathbb{F}_q)$ is fixed; or
 \begin{equation}\label{maximumcliqueg}
 \mbox{${\cal M}=\left(\mathbb{F}_{q}^{m\times (d-1)}, 0 \right)Q:=\left\{(Y,0)Q : Y\in\mathbb{F}_{q}^{m\times (d-1)} \right\}$}
 \end{equation}
 with $n=m$, where $Q\in GL_m(\mathbb{F}_q)$ is fixed.
 \end{lem}
\proof
Let $\varphi$ be the graph isomorphism (\ref{isomorphismattenuated}) from the adjacency graph of $\left( \mathcal{A}_{m-1}, \mathcal{A}_{m}, \subseteq \right)$ to
the bilinear forms graph $\Gamma_2(\mathbb{F}_{q}^{m\times n})$.
For any $A,B\in\mathbb{F}_{q}^{m\times n}$, we have ${\rm rank}(A-B)<d$ (where $2\leq d\leq m$) if and only if ${\rm dim}(\varphi^{-1}(A)\cap\varphi^{-1}(B))\geq m-d+1$.
Thus, ${\cal M}$ is a maximum clique of $\Gamma_d(\mathbb{F}_{q}^{m\times n})$ if and only if
$\varphi^{-1}({\cal M})$ has the  property that ${\rm dim}(\gamma \cap \delta)\geq m-d+1$ for all $\gamma, \delta$ in $\varphi^{-1}({\cal M})$ and
 $|\varphi^{-1}({\cal M})|={\rm max}\{|\mathscr{F}|\}$ for all collection $\mathscr{F}$ of elements of
$\mathcal{A}_{m}$ with the property that ${\rm dim}(\gamma \cap \delta)\geq m-d+1$ for all $\gamma, \delta$ in $\mathscr{F}$.

Let ${\cal M}$ be a maximum clique  containing $0$ in $\Gamma_d(\mathbb{F}_{q}^{m\times n})$. Write $W=(I_n,0)$. By Lemma \ref{Rectangular-PID2-4a} and above result,
$\varphi^{-1}({\cal M})$ is of the form either
{\rm (a)}   $\varphi^{-1}({\cal M})$  consists of all elements of $\mathcal{A}_{m}$ which contain a fixed $(m-d+1)$-dimensional subspace $U$ with $U\cap W=0$, or
 {\rm (b)} $m=n$ and  $\varphi^{-1}({\cal M})$  is the set of all elements of  $\mathcal{A}_{m}$ contained in a fixed  $(m+d-1)$-dimensional subspace $U'$ with
 ${\rm dim}(U'\cap W)=d-1$. In the case (a), it is easy to see that  ${\cal M}$  is of  the  form  (\ref{maxiumclique35trgr}). Now, we assume the case (b) happens.
Since ${\rm dim}(U'\cap W)=d-1$, by appropriate elementary operations of matrix, we may assume with no loss of generality that $\small U'=\left(
                                                                                                  \begin{array}{ccc}
                                                                                                    I_{d-1} &0& 0 \\
                                                                                                    0 &  0&I_{m} \\
                                                                                                  \end{array}
                                                                                                \right)$.
Then $\varphi^{-1}({\cal M})=\left\{(Y,0,I_m): Y\in\mathbb{F}_{q}^{m\times (d-1)} \right\}$. Thus (\ref{maximumcliqueg}) holds.
$\qed$

However, the proof of Lemma \ref{Rectangular-PID2-4} cannot be generalized to the case of $\mathbb{Z}_{p^s}$. In order to generalize  Lemma \ref{Rectangular-PID2-4}
to the case of $\mathbb{Z}_{p^s}$, we need a new method.

\begin{lem}\label{rankinvertible935tet} \
Let $k\geq 2$ and $A\in\mathbb{Z}_{p}^{k\times k}$. Then there exists $B\in GL_k(\mathbb{Z}_{p})$ such that $A-B$ is invertible.
\end{lem}
\proof
Without losing generality, we may assume that $A=I_k$.
When $k$ is an even number, we have $\small I_k-\left(
                                           \begin{array}{cc}
                                             I_{k/2} & I_{k/2} \\
                                             I_{k/2} & 0 \\
                                           \end{array}
                                         \right)$ is invertible, and hence this lemma holds.
From now on we assume that  $k$ is an odd number. When $k=3$, there is $\scriptsize B_1:=\left(
    \begin{array}{ccc}
      1 & 1 & 0\\
      1& 0 & 1 \\
      0 & 1 & 0 \\
    \end{array}
  \right)\in GL_3(\mathbb{Z}_{p})$ such that $I_3-B_1$ is invertible.
 When $k=2r+1$ ($r\geq 2$), by the case of even, there is $B_2\in GL_{2r-2}(\mathbb{Z}_{p})$ such that $I_{2r-2}-B_1$ is invertible,
 and hence $I_{2r+1}-{\rm diag}(B_1,B_2)$ is invertible.
$\qed$

\begin{thm}\label{maximumclique-ggggbbb} \  Let $n\geq m$, and let ${\cal M}$ be a maximum clique of $\Gamma_d(\mathbb{Z}_{p^s}^{m\times n})$.
Then  ${\cal M}$ is of the form either
 \begin{equation}\label{maximumcliquete408}
 \mbox{${\cal M}=P\left(\begin{array}{c}
                              \mathbb{Z}_{p^s}^{(d-1)\times n} \\
                              0 \\
                            \end{array}
                          \right)+B:=\left\{P\left(
                            \begin{array}{c}
                              X\\
                              0 \\
                            \end{array}
                          \right)+B: X\in\mathbb{Z}_{p^s}^{(d-1)\times n} \right\}$,}
 \end{equation}
  where $P\in GL_m(\mathbb{Z}_{p^s})$ and $B\in\mathbb{Z}_{p^s}^{m\times n}$ are fixed; or
 \begin{equation}\label{maximumcliq786dfdsj0gf}
 \mbox{${\cal M}=\left(\mathbb{Z}_{p^s}^{m\times (d-1)}, 0 \right)Q+B:=\left\{(Y,0)Q+B : Y\in\mathbb{Z}_{p^s}^{m\times (d-1)} \right\}$}
 \end{equation}
 with $n=m$, where $Q\in GL_m(\mathbb{Z}_{p^s})$ and $B\in\mathbb{Z}_{p^s}^{m\times m}$ are fixed.
\end{thm}
\proof
When $s=1$, this theorem holds by Lemma \ref{Rectangular-PID2-4}.
By \cite[Theorem 3.6]{Huang-2017}, this theorem holds if $d=2$. Thus, from now on we assume that $s\geq 2$ and $m\geq d\geq 3$.
We prove this theorem by induction on $s$. Assume that this theorem holds for $s-1$. We prove that it holds for $s$ as follows.

Let  ${\cal M}=\{A_1,\ldots, A_{\omega}\}$ be a maximum clique of $\Gamma_d(\mathbb{Z}_{p^s}^{m\times n})$.  By (\ref{biline13nvhiu6})
we have $\omega=|{\cal M}|=p^{sn(d-1)}$. By
the bijection $X\mapsto X-B$, we can assume that ${\cal M}$ contains $0$ (i.e., $B=0$).

Let  $\pi: \mathbb{Z}_{p^s}^{m\times n}\rightarrow \mathbb{Z}_{p}^{m\times n}$ be the natural surjection (\ref{naturalsurjection}).
For any vertex $X$ in $\pi(\mathcal{M})$, let $\pi^{-1}(X)$ denote the  preimages of $X$ in $\mathcal{M}$, i.e.,
$\pi^{-1}(X)=\left\{Y\in \mathcal{M}: \pi(Y)=X \right\}$.
Suppose $\{\pi(A_{i_1}), \ldots, \pi(A_{i_h})\}$ is the set of all different elements in $\{\pi(A_1), \ldots, \pi(A_\omega)\}$.
Then   $\pi({\cal M})=\{\pi(A_{i_1}), \ldots, \pi(A_{i_h})\}$ is a clique of $\Gamma_d(\mathbb{Z}_{p}^{m\times n})$.
It follows from (\ref{geneytr675vxcwq128}) that
\begin{equation}\label{gd54697qekh9c}
h\leq \omega\left(\Gamma_d(\mathbb{Z}_{p}^{m\times n})\right)=p^{n(d-1)}.
\end{equation}
 Clearly, ${\cal M}$ has a  partition into $h$  cliques
$$
{\cal M}={\cal M}_1\cup {\cal M}_2\cup \cdots \cup {\cal M}_h,
$$
 where $\pi({\cal M}_t)=\{\pi(A_{i_t})\}$, $t=1,\ldots,h$.
Thus
\begin{equation}\label{bvcqe235jgh75df}
|{\cal M}|=p^{sn(d-1)}=\sum_{t=1}^h|M_t|.
\end{equation}

Let $n_t=|{\cal M}_t|$ and let
$${\cal M}_t=\left\{\pi(A_{i_t})+B_{1t}, \pi(A_{i_t})+B_{2t}, \ldots, \pi(A_{i_t})+B_{n_t,t}\right\},$$
 where $\left\{B_{1t}, B_{2t}, \ldots, B_{n_t,t}\right\}$ is a clique in $J_{p^s}^{m\times n}$,
$t=1,\ldots, h$. Write $B_{jt}=pC_{jt}$ where $C_{jt}\in\mathbb{Z}_{p^{s-1}}^{m\times n}$, $j=1,\ldots,n_t$, $t=1,\ldots,h$.
By Lemma \ref{64gdgd7wrww}, $\left\{C_{1t}, C_{2t}, \ldots, C_{n_t,t}\right\}$ is a clique of $\Gamma_d(\mathbb{Z}_{p^{s-1}}^{m\times n})$,
 $t=1,\ldots,h$. By (\ref{biline13nvhiu6}), we get
\begin{equation}\label{gf5nvj979z0hjgy}
n_t=|{\cal M}_t|\leq \omega\left(\Gamma_d(\mathbb{Z}_{p^{s-1}}^{m\times n})\right)=p^{(s-1)n(d-1)},  \ \ t=1,\ldots,h.
\end{equation}
Thus $p^{sn(d-1)}=|{\cal M}|\leq hp^{(s-1)n(d-1)}$, and hence $p^{n(d-1)}\leq h$. By (\ref{gd54697qekh9c}), we obtain
\begin{equation}\label{GFDWE546MBCXSF7}
h=p^{n(d-1)}=\omega(\mathbb{Z}_{p}^{m\times n}).
\end{equation}
Therefore, $\pi({\cal M})$ is a maximum clique containing $0$ in $\Gamma_d(\mathbb{Z}_{p}^{m\times n})$.

By  Lemma \ref{Rectangular-PID2-4},  we have either
$\small\pi({\cal M})=P_0\left(
\begin{array}{c}
\mathbb{Z}_{p}^{(d-1)\times n} \\
0 \\
\end{array}\right)$, or  $\pi({\cal M})=\left(\mathbb{Z}_p^{m\times (d-1)}, 0 \right)Q_0$ with $n=m$, where $P_0,Q_0$ are invertible matrices over $\mathbb{Z}_{p}$.
Thus, ${\cal M}$ contains a matrix $E$ of inner rank $d-1$, such that either $\pi(E)=P_0{\rm diag}(I_{d-1},0)$ or $\pi(E)={\rm diag}(I_{d-1},0)Q_0$.
 Using appropriate elementary operations of matrix, without losing generality, we may assume that
\begin{equation}\label{gfd54dnca68hldx0}
E={\rm diag}(I_{d-1},0)=\pi(E)\in \mathcal{M}.
\end{equation}
Since $\scriptsize E=P_0\left(
\begin{array}{c}
X_1 \\
0 \\
\end{array}\right)$ where $X_1\in\mathbb{Z}_{p}^{(d-1)\times n}$, or  $E=\left(Y_1, 0 \right)Q_0$ with $n=m$, where $Y_1\in \mathbb{Z}_p^{m\times (d-1)}$,
it follows  that $\scriptsize P_0=\left(
\begin{array}{cc}
 P_{11} &  * \\
    0 & * \\
    \end{array}
 \right)$ or $\scriptsize Q_0=\left(
\begin{array}{cc}
 Q_{11} & 0 \\
    * & * \\
    \end{array}
 \right)$, where $P_{11}, Q_{11}\in GL_{d-1}(\mathbb{Z}_{p})$.
Thus,  $\pi({\cal M})$ is of the form either
\begin{equation}\label{ghf435nvz9078rf}
\pi({\cal M})=\left(
\begin{array}{c}
\mathbb{Z}_{p}^{(d-1)\times n} \\
0 \\
\end{array}\right),
\end{equation}
or
\begin{equation}\label{675fd425klhlcvx07}
\mbox{$\pi({\cal M})=\left(\mathbb{Z}_p^{m\times (d-1)}, 0 \right)$ \ with \  $n=m$.}
\end{equation}

By (\ref{bvcqe235jgh75df})--(\ref{GFDWE546MBCXSF7}), it is easy to see that
$$
|{\cal M}_t|=n_t=\omega\left(\Gamma_d(\mathbb{Z}_{p^{s-1}}^{m\times n})\right)=p^{(s-1)n(d-1)}, \ \ t=1,\ldots, h.
$$
 Thus
\begin{equation}\label{maxicliq5436uvbvl087}
\left|\pi^{-1}(\pi(A_{i_t}))\right|=|{\cal M}_t|=p^{(s-1)n(d-1)}=\omega\left(\Gamma_d(\mathbb{Z}_{p^{s-1}}^{m\times n})\right),  \ t=1,\ldots, h.
\end{equation}
It follows that  $\left\{C_{1t}, C_{2t}, \ldots, C_{n_t,t}\right\}$ is a maximum  clique of $\Gamma_d(\mathbb{Z}_{p^{s-1}}^{m\times n})$, $t=1,\ldots, h$.
In other words, there exists a maximum  clique $\mathcal{C}_t$ of $\Gamma_d(\mathbb{Z}_{p^{s-1}}^{m\times n})$  such that
\begin{equation}\label{maxicliq5436uvbvl088}
\pi^{-1}(\pi(A_{i_t}))={\cal M}_t=\pi(A_{i_t})+\mathcal{C}_tp,  \ t=1,\ldots, h.
\end{equation}
By the induction hypothesis, $\mathcal{C}_t$ is of the form either (\ref{maximumcliquete408}), or (\ref{maximumcliq786dfdsj0gf}) with $n=m$.
Thus,
\begin{equation}\label{maxicliq5436uvbvl090}
\mbox{$\mathcal{C}_t=\left(
\begin{array}{c}
 P_{t1}\mathbb{Z}_{p^{s-1}}^{(d-1)\times n}p  \\
 P_{t2}\mathbb{Z}_{p^{s-1}}^{(d-1)\times n}p \\
 \end{array}
 \right)+B_t$, \ or \ $\mathcal{C}_t=\left(\mathbb{Z}_{p^{s-1}}^{m\times (d-1)}p Q_{t1},  \mathbb{Z}_{p^{s-1}}^{m\times (n-d+1)}p Q_{t2}\right)+B_t$,}
\end{equation}
where $P_{t1}, Q_{t1}\in\mathbb{Z}_{p^{s-1}}^{(d-1)\times (d-1)}$,  $\scriptsize\left(
\begin{array}{c}
P_{t1}  \\
 P_{t2}  \\
  \end{array}
  \right)\in \mathbb{Z}_{p^{s-1}}^{m\times (d-1)}$ has a left inverse,  $(Q_{t1}, Q_{t2})\in\mathbb{Z}_{p^{s-1}}^{(d-1)\times n}$
  has a right inverse, and $B_t\in\mathbb{Z}_{p^{s-1}}^{m\times n}$, $t=1,\ldots, h$.

Let $A_{i_1}=0$. Then $\pi(0)=0$, $\pi^{-1}(0)$ is of the form (\ref{maxicliq5436uvbvl088}) and  $\mathcal{C}_1$ contains $0$.
  By  (\ref{maxicliq5436uvbvl090}), we have either
\begin{equation}\label{FGD543JKCX965G}
    \pi^{-1}(0)= \left\{\left(
     \begin{array}{cc}
    P_{11}Xp & P_{11}Yp  \\
     P_{12}Xp  & P_{12}Yp  \\
     \end{array}
    \right): X\in\mathbb{Z}_{p^{s-1}}^{(d-1)\times (d-1)}, Y\in\mathbb{Z}_{p^{s-1}}^{(d-1)\times (n-d+1)}
   \right\},
\end{equation}
or
\begin{equation}\label{86gdgrevbcml3JKCX9G}
\mbox{$\pi^{-1}(0)= \left\{\left(
\begin{array}{cc}
 XpQ_{11} & XpQ_{12}  \\
  Wp Q_{11}  & WpQ_{12}  \\
  \end{array}
 \right): X\in\mathbb{Z}_{p^{s-1}}^{(d-1)\times (d-1)}, W\in\mathbb{Z}_{p^{s-1}}^{(n-d+1)\times (d-1)}
   \right\}$ \ \ with \ $n=m$.}
 \end{equation}

Suppose  (\ref{FGD543JKCX965G}) holds.
Since (\ref{gfd54dnca68hldx0}) and  $\small\left(\begin{array}{cc}
 0 & P_{11}Yp  \\
 0  & P_{12}Yp  \\
 \end{array}\right)\sim \left(
                                \begin{array}{cc}
                                  I_{d-1} & 0 \\
                                  0  & 0  \\
                                \end{array}
                              \right)$,
$\small\rho \left(
\begin{array}{cc}
-I_{d-1} & P_{11}Yp  \\
0 & P_{12}Yp  \\
\end{array}\right)\leq d-1$ for all  $ Y\in\mathbb{Z}_{p^{s-1}}^{(d-1)\times (n-d+1)}$. By (\ref{rank-0008}), we get that $P_{12}=0$ and $P_{11}$ is invertible.
Therefore
\begin{equation}\label{maximlif425rwrsd02}
\pi^{-1}(0)= \left\{\left(\begin{array}{cc}
Xp & Yp  \\
 0  & 0  \\
  \end{array}
  \right): (X,Y)\in\mathbb{Z}_{p^{s-1}}^{(d-1)\times n}\right\}.
\end{equation}

Suppose  (\ref{86gdgrevbcml3JKCX9G}) holds.  We have similarly that
\begin{equation}\label{hfg65ryuvxcb0akk5}
\mbox{$\pi^{-1}(0)= \left\{\left(\begin{array}{cc}
Xp & 0  \\
Wp  & 0  \\
  \end{array}
  \right): \left(\begin{array}{c}
X  \\
W  \\
  \end{array}
  \right)\in\mathbb{Z}_{p^{s-1}}^{m\times (d-1)}\right\}$ \ \ with \ $n=m$.}.
\end{equation}
Thus, $\pi^{-1}(0)$ is of the form either (\ref{maximlif425rwrsd02}) or (\ref{hfg65ryuvxcb0akk5}) with $n=m$.

We distinguish the following two cases to prove this theorem.

{\bf Case 1.} \  $\pi({\cal M})$ is of the form (\ref{ghf435nvz9078rf}).

 First, we show that (\ref{maximlif425rwrsd02}) holds. Let $\scriptsize D=\left(\begin{array}{cc}
0 &  I_{d-1}  \\
 0  & 0  \\
  \end{array}
  \right)\in \pi({\cal M})$.  By (\ref{maxicliq5436uvbvl088}), there is a maximum  clique $\mathcal{C}_D$
of $\Gamma_d(\mathbb{Z}_{p^{s-1}}^{m\times n})$  such that
$$\pi^{-1}\left(\begin{array}{cc}
0 & I_{d-1}  \\
 0  & 0  \\
  \end{array}
  \right)=
\left\{\left(\begin{array}{cc}
 X_1p & I_{d-1} +X_2p \\
 X_3p & X_4p  \\
  \end{array}
  \right): \left(
 \begin{array}{cc}
   X_1 & X_2 \\
     X_3 & X_4 \\
      \end{array}
   \right)\in \mathcal{C}_D
   \right\}.$$
Suppose that $\pi^{-1}(0)$ is of the form (\ref{hfg65ryuvxcb0akk5}) with $n=m$. Then
$\small\left(\begin{array}{cc}
 X_1p & I_{d-1} +X_2p \\
 X_3p & X_4p  \\
  \end{array}
  \right)\sim \left(\begin{array}{cc}
Xp & 0  \\
Wp  & 0  \\
  \end{array}
  \right)$ for all $X\in\mathbb{Z}_{p^{s-1}}^{(d-1)\times (d-1)}$ and $W\in\mathbb{Z}_{p^{s-1}}^{(m-d+1)\times (d-1)}$.
Thus, we can choose $X$ and $W$ such that
$$\rho\left(\begin{array}{ccc}
 0&* & I_{d-1} +X_2p \\
 W_1p &*& X_4p  \\
  \end{array}
  \right)=\rho\left(\begin{array}{ccc}
 0&0& I_{d-1} +X_2p \\
 W_1p &*& 0 \\
  \end{array}
  \right)\leq d-1,
  $$
where $0\neq W_1\in\mathbb{Z}_{p^{s-1}}^{(m-d+1)\times 1}$. By (\ref{rank-0008}), this is a contradiction.
Thus  $\pi^{-1}(0)$ must be of the form (\ref{maximlif425rwrsd02}).

Let $e_i$ be the $i$-th  column  of $I_{d-1}$.  By (\ref{maxicliq5436uvbvl088}) and (\ref{maxicliq5436uvbvl090}), there is a maximum  clique $\mathcal{C}_{e_i}$
of $\Gamma_d(\mathbb{Z}_{p^{s-1}}^{m\times n})$  such that
$$\pi^{-1}\left(\begin{array}{cc}
0 &  e_1  \\
 0  & 0  \\
  \end{array}
  \right)=
\left\{\left(\begin{array}{cc}
 X_1p & e_1 +Y_3p \\
 X_2p & Y_4p  \\
  \end{array}
  \right): \left(
 \begin{array}{cc}
   X_1 & Y_3  \\
     X_2 & Y_4 \\
      \end{array}
   \right)\in \mathcal{C}_{e_1}
   \right\}.$$
Let $T$ be any $(n-1)\times (n-1)$ permutation matrix over $\mathbb{Z}_p$. Then $\rho\left({\rm diag}(0,I_{d-2})Tp, \ e_1\right)=d-1$.
Write $X_2=X_2'T$, $X_2\in\mathbb{Z}_{p^{s-1}}^{(m-d+1)\times (n-1)}$.
Using (\ref{maximlif425rwrsd02}), we get
$$\left(\begin{array}{cc}
 X_1p & e_1+Y_3p \\
  X_2p & Y_4p \\
  \end{array}
  \right)\sim \left(\begin{array}{cc}
 X_1p-{\rm diag}(0,I_{d-2})Tp & Y_3p \\
  0  & 0 \\
  \end{array}
  \right).$$
Thus
$$\rho\left(\begin{array}{cc}
{\rm diag}(0,I_{d-2})Tp & e_1  \\
 X_2'Tp & Y_4p   \\
  \end{array}
  \right)=\rho\left(\begin{array}{cc}
{\rm diag}(0,I_{d-2})Tp & e_1  \\
(X'_{21},0 )Tp & 0  \\
  \end{array}
  \right)\leq d-1,
$$
where $X_{21}'\in\mathbb{Z}_{p^{s-1}}^{(m-d+1)\times (n-d+1)}$. It follows from (\ref{rank-0008}) that $X_{21}'=0$, and hence the matrix
$X_2$ has $n-d+1$ columns to be zeros. Since the permutation matrix $T$ is arbitrary, every column of $X_2$ must be zero, and hence $X_2=0$. Then
$$\pi^{-1}\left(\begin{array}{cc}
0 &  e_1  \\
 0 & 0\\
  \end{array}
  \right)=
\left\{\left(\begin{array}{cc}
 X_1p & e_1 +Y_3p \\
 0  & Y_4p\\
  \end{array}
  \right): \left(
 \begin{array}{cc}
   X_1 & Y_3  \\
   0 & Y_4 \\
      \end{array}
   \right)\in \mathcal{C}_{e_1}
   \right\}.$$
Since $\small\left(\begin{array}{cc}
X_1p & e_1+Y_3p  \\
 0 & Y_4p  \\
  \end{array}
  \right) \sim
\left(\begin{array}{cc}
I_{d-1} &  0  \\
 0  & 0  \\
  \end{array}
  \right)$, it is easy to see that $Y_4p=0$. Therefore, we obtain that
\begin{equation}\label{sfd425badai9b}
\pi^{-1}\left(\begin{array}{cc}
 0& e_1   \\
 0 & 0  \\
  \end{array}
  \right)=
\left\{\left(\begin{array}{cc}
 Xp & e_1 +Yp  \\
 0 & 0  \\
  \end{array}
  \right): X\in\mathbb{Z}_{p^{s-1}}^{(d-1)\times (n-1)},  Y\in\mathbb{Z}_{p^{s-1}}^{(d-1)\times 1}
   \right\}.
\end{equation}

Similarly, we can prove that
\begin{equation}\label{sfd425badai9iii}
\pi^{-1}\left(\begin{array}{cc}
 0& e_i   \\
 0 & 0  \\
  \end{array}
  \right)=\left\{\left(\begin{array}{cc}
 Xp & e_i +Yp  \\
 0 & 0  \\
  \end{array}
  \right): X\in\mathbb{Z}_{p^{s-1}}^{(d-1)\times (n-1)},  Y\in\mathbb{Z}_{p^{s-1}}^{(d-1)\times 1}\right\}, \ \  i=1,\ldots, d-1.
\end{equation}

Let $A\in GL_{d-1}(\mathbb{Z}_{p})$.  By (\ref{maxicliq5436uvbvl088}), there is a maximum  clique $\mathcal{C}_A$
of $\Gamma_d(\mathbb{Z}_{p^{s-1}}^{m\times n})$  such that
$$\mbox{$\pi^{-1}\left(\begin{array}{ccc}
A & 0  \\
 0 & 0  \\
  \end{array}
  \right)=
\left\{\bordermatrix{
&_{d-1} & _{n-d} & _{1} \cr
& A+ X_1p & Y_1p& Z_1p\cr
& X_2p  & Y_2p& Z_2p  \cr }
: \left(
                                \begin{array}{ccc}
                                  X_1& Y_1& Z_1  \\
                                  X_2  & Y_2& Z_2 \\
                                \end{array}
                              \right)\in \mathcal{C}_A
   \right\}.$}$$
By (\ref{maximlif425rwrsd02}), we have
$\mbox{\small$\bordermatrix{
&_{d-1} & _{n-d} & _{1} \cr
& A+ X_1p & Y_1p& Z_1p\cr
& X_2p  & Y_2p& Z_2p  \cr } \ \sim
\bordermatrix{
&_{d-1} & _{n-d} & _{1} \cr
&  X_1p & Y_1p& Z_1p\cr
& 0  & 0& 0 \cr }$}$.  Thus
 $\small\rho\left(\begin{array}{ccc}
 A & 0& 0\\
 X_2p  & Y_2p& Z_2p  \\
  \end{array}
  \right)\leq d-1$, and hence  $Y_2p=0$ and $Z_2p=0$.
 On the other hand, from  (\ref{sfd425badai9iii}) we get
 $$\mbox {$
\bordermatrix{
&_{d-1} & _{n-d} & _{1} \cr
& A+ X_1p & Y_1p& Z_1p\cr
& X_2p  & 0& 0  \cr } \ \sim
\bordermatrix{
&_{d-1} & _{n-d} & _{1} \cr
&  X_1p & Y_1p& e_i+Z_1p\cr
& 0  & 0& 0 \cr }$, \ $i=1, \ldots, d-1$.}$$
Consequently
 $\small\rho\left(\begin{array}{ccc}
 A & 0& -e_i\\
 X_2p  & 0 & 0  \\
  \end{array}
  \right)=\rho\left(\begin{array}{ccc}
 A & 0& -e_i\\
  0 & 0 & X_2A^{-1}e_ip  \\
  \end{array}
  \right)\leq d-1$, and hence  $X_2A^{-1}e_i=0$ by (\ref{rank-0008}), $i=1, \ldots, d-1$. Therefore $X_2=0$.
  Then we obtain
\begin{equation}\label{s43REd425baiWEQR}
\mbox{$\pi^{-1}\left(\begin{array}{cc}
 A& 0   \\
 0 & 0  \\
  \end{array}
  \right)=
\left\{\left(\begin{array}{cc}
 A+ Xp &Yp  \\
 0 & 0  \\
  \end{array}
  \right): X\in\mathbb{Z}_{p^{s-1}}^{(d-1)\times (d-1)}, Y\in\mathbb{Z}_{p^{s-1}}^{(d-1)\times (n-d+1)}
   \right\}$, \ \ $A\in GL_{d-1}(\mathbb{Z}_{p})$.}
\end{equation}

Let $A\in GL_{d-1}(\mathbb{Z}_{p})$.  Write $F_i=(A,e_i)$. By (\ref{maxicliq5436uvbvl088}), there is a maximum  clique $\mathcal{C}_{F_i}$
of $\Gamma_d(\mathbb{Z}_{p^{s-1}}^{m\times n})$  such that
$$\mbox{$\pi^{-1}\left(\begin{array}{ccc}
A & 0 &e_i \\
 0 & 0& 0  \\
  \end{array}
  \right)=
\left\{\bordermatrix{
&_{d-1} & _{n-d} & _{1} \cr
& A+ X_1p & Y_1p& e_i+Z_1p\cr
& X_2p  & Y_2p& Z_2p  \cr }
: \left(
                                \begin{array}{ccc}
                                  X_1& Y_1& Z_1  \\
                                  X_2  & Y_2& Z_2 \\
                                \end{array}
                              \right)\in \mathcal{C}_{F_i}
   \right\}.$}$$
From (\ref{sfd425badai9iii}) we have $\mbox{\small$\bordermatrix{
&_{d-1} & _{n-d} & _{1} \cr
& A+ X_1p & Y_1p& e_i+Z_1p\cr
& X_2p  & Y_2p& Z_2p  \cr } \ \sim \bordermatrix{
&_{d-1} & _{n-d} & _{1} \cr
& X_1p & Y_1p& e_i+Z_1p\cr
& 0 & 0& 0\cr }$}$, which implies that $Y_2p=0$ and $Z_2p=0$. By (\ref{maximlif425rwrsd02}), we get
$\mbox{\small$\bordermatrix{
&_{d-1} & _{n-d} & _{1} \cr
& A+ X_1p & Y_1p& e_i+Z_1p\cr
& X_2p  & 0& 0 \cr } \ \sim \bordermatrix{
&_{d-1} & _{n-d} & _{1} \cr
& X_1p & Y_1p& Z_1p\cr
& 0 & 0& 0\cr }$}$, and thus
$$\rho\left(
    \begin{array}{ccc}
      A & 0& e_i \\
    X_2p  & 0& 0 \\
    \end{array}
  \right)=\rho\left(
    \begin{array}{ccc}
      A & 0& e_i \\
     0 & 0& -X_2A^{-1}e_ip  \\
    \end{array}
  \right)\leq d-1.$$
By (\ref{rank-0008}), one has $X_2A^{-1}e_i=0$. Using (\ref{sfd425badai9iii}) again,  we get
 $$\mbox{$\bordermatrix{
&_{d-1} & _{n-d} & _{1} \cr
& A+ X_1p & Y_1p& e_i+Z_1p\cr
& X_2p  & Y_2p& Z_2p  \cr } \ \sim \bordermatrix{
&_{d-1} & _{n-d} & _{1} \cr
& X_1p & Y_1p& e_j+Z_1p\cr
& 0 & 0& 0\cr }$}, \ \ j=1,\ldots,d-1, $$
and hence
$$\rho\left(
    \begin{array}{ccc}
      A & 0& e_i-e_j \\
    X_2p  & 0& 0 \\
    \end{array}
  \right)=\rho\left(
    \begin{array}{ccc}
      A & 0& e_i-e_j \\
     0 & 0& -X_2A^{-1}(e_i-e_j)p  \\
    \end{array}
  \right)\leq d-1. \ \ j=1,\ldots,d-1.$$
It follows from $X_2A^{-1}e_i=0$ that $X_2A^{-1}e_j=0$, $j=1,\ldots,d-1$, which implies that $X_2=0$.
Thus, for any $A\in GL_{d-1}(\mathbb{Z}_{p})$, we have
\begin{equation}\label{65GDHF8r6vbs42C}
\mbox{\small$\pi^{-1}\left(\begin{array}{ccc}
A & 0 &e_i \\
 0 & 0& 0  \\
  \end{array}
  \right)=
\left\{\bordermatrix{
&_{d-1} & _{n-d} & _{1} \cr
& A+X_1p & Y_1p& e_i+Z_1p\cr
& 0 & 0& 0 \cr }
: (X_1, Y_1,Z_1)\in\mathbb{Z}_{p^{s-1}}^{(d-1)\times n}
   \right\}$, \ \  $i=1,\ldots, d-1$,}
\end{equation}

Let $B\in\mathbb{Z}_p^{(d-1)\times(d-1)}$.  By (\ref{maxicliq5436uvbvl088}), there is a maximum  clique $\mathcal{C}_B$
of $\Gamma_d(\mathbb{Z}_{p^{s-1}}^{m\times n})$  such that
$$\mbox{$\pi^{-1}\left(\begin{array}{ccc}
B & 0  \\
 0 & 0  \\
  \end{array}
  \right)=
\left\{\bordermatrix{
&_{d-1} & _{n-d} & _{1} \cr
& B+ X_1p & Y_1p& Z_1p\cr
& X_2p  & Y_2p& Z_2p  \cr }
: \left(
                                \begin{array}{ccc}
                                  X_1& Y_1& Z_1  \\
                                  X_2  & Y_2& Z_2 \\
                                \end{array}
                              \right)\in \mathcal{C}_B
   \right\}.$}$$
By Lemma \ref{rankinvertible935tet}, exists $A\in GL_{d-1}(\mathbb{Z}_{p})$ such that $B-A$ is invertible. Applying (\ref{s43REd425baiWEQR}), we get
$$\mbox{$\bordermatrix{
&_{d-1} & _{n-d} & _{1} \cr
& B+ X_1p & Y_1p& Z_1p\cr
& X_2p  & Y_2p& Z_2p  \cr } \ \sim \bordermatrix{
&_{d-1} & _{n-d} & _{1} \cr
& A+X_1p & Y_1p& Z_1p\cr
& 0 & 0& 0\cr }$},$$
 which implies that $Y_2p=0$ and $Z_2p=0$. On the other hand, from (\ref{65GDHF8r6vbs42C}) we have
$$\mbox{$\bordermatrix{
&_{d-1} & _{n-d} & _{1} \cr
& B+ X_1p & Y_1p& Z_1p\cr
& X_2p  & 0& 0 \cr } \ \sim \bordermatrix{
&_{d-1} & _{n-d} & _{1} \cr
& A+X_1p & Y_1p& e_i+Z_1p\cr
& 0 & 0& 0\cr }$},$$
it follows that
$$\rho\left(
    \begin{array}{ccc}
     B-A & 0& -e_i \\
    X_2p  & 0& 0 \\
    \end{array}
  \right)=\rho\left(
    \begin{array}{ccc}
      B-A & 0& -e_i \\
     0 & 0& X_2(B-A)^{-1}e_ip  \\
    \end{array}
  \right)\leq d-1, \ \ i=1,\ldots, d-1.$$
By (\ref{rank-0008}), we obtain $X_2(B-A)^{-1}=0$ and hence $X_2=0$. Then we have proved that
\begin{equation}\label{f21jgxad69gnklkbb}
\pi^{-1}\left(\begin{array}{cc}
B & 0  \\
 0 & 0  \\
  \end{array}\right)=
 \left\{\left(\begin{array}{cc}
B+ Xp & Yp  \\
 0 & 0  \\
  \end{array}
  \right): X\in\mathbb{Z}_{p^{s-1}}^{(d-1)\times (d-1)},  Y\in\mathbb{Z}_{p^{s-1}}^{(d-1)\times (n-d+1)} \right\}, \ \ B\in\mathbb{Z}_p^{(d-1)\times(d-1)}.
\end{equation}
In particular, we have
\begin{equation}\label{gfd54bcgsi9i4}
\mbox{$\pi^{-1}\left(\begin{array}{cc}
 e_i & 0  \\
 0 & 0  \\
  \end{array}
  \right)=
\left\{\left(\begin{array}{cc}
 e_i+ Xp & Yp  \\
 0 & 0  \\
  \end{array}
  \right): X\in\mathbb{Z}_{p^{s-1}}^{(d-1)\times 1}, Y\in\mathbb{Z}_{p^{s-1}}^{(d-1)\times (n-1)}
   \right\}$, \ \  $i=1,\ldots, d-1$.}
\end{equation}

Using (\ref{gfd54bcgsi9i4}) and (\ref{maximlif425rwrsd02}), similar to the proof of (\ref{s43REd425baiWEQR}), we can get that
\begin{equation}\label{iWE543ghd65Qi8}
\mbox{$\pi^{-1}\left(\begin{array}{cc}
 0& A  \\
 0 & 0  \\
  \end{array}
  \right)=\left\{\left(\begin{array}{cc}
 Xp &A+Yp \\
 0 & 0  \\
  \end{array}
  \right): X\in\mathbb{Z}_{p^{s-1}}^{(d-1)\times (n-d+1)}, Y\in\mathbb{Z}_{p^{s-1}}^{(d-1)\times (d-1)}
   \right\}$, \ \ $A\in GL_{d-1}(\mathbb{Z}_{p})$.}
\end{equation}

Let $A\in GL_{d-1}(\mathbb{Z}_{p})$.  Similar to the proof of (\ref{65GDHF8r6vbs42C}), by  (\ref{maximlif425rwrsd02}) and (\ref{gfd54bcgsi9i4}) we can prove that
\begin{equation}\label{ytgwqlbnbc85f}
\mbox{\small$\pi^{-1}\left(\begin{array}{ccc}
 e_i& 0 &A \\
 0 & 0& 0  \\
  \end{array}
  \right)=
\left\{\bordermatrix{
&_{1} & _{n-d} & _{d-1} \cr
& e_i+X_1p & Y_1p& A+Z_1p\cr
& 0 & 0& 0 \cr }
: (X_1, Y_1,Z_1)\in\mathbb{Z}_{p^{s-1}}^{(d-1)\times n}
   \right\}$, \ \  $i=1,\ldots, d-1$.}
\end{equation}

Using (\ref{iWE543ghd65Qi8}) and (\ref{ytgwqlbnbc85f}), similar to the proof of (\ref{f21jgxad69gnklkbb}), we have
\begin{equation}\label{5HF713SFBlk3}
\pi^{-1}\left(\begin{array}{cc}
0 & B  \\
 0 & 0  \\
  \end{array}
  \right)=
\left\{\left(\begin{array}{cc}
Xp & B+Yp \\
 0 & 0  \\
  \end{array}
  \right): X\in\mathbb{Z}_{p^{s-1}}^{(d-1)\times (n-d+1)}, Y\in\mathbb{Z}_{p^{s-1}}^{(d-1)\times (d-1)}\right\},  \ \ B\in\mathbb{Z}_p^{(d-1)\times(d-1)}.
\end{equation}

Now, let $A\in GL_{d-1}(\mathbb{Z}_{p})$ and $B_1\in\mathbb{Z}_p^{(d-1)\times(n-d+1)}$. Write $H=(A,B_1)$.
 By (\ref{maxicliq5436uvbvl088}), there is a maximum  clique $\mathcal{C}_{H}$
of $\Gamma_d(\mathbb{Z}_{p^{s-1}}^{m\times n})$  such that
$$\mbox{$\pi^{-1}\left(\begin{array}{cc}
A  & B_1 \\
0& 0  \\
  \end{array}
  \right)=
\left\{\left(\begin{array}{cc}
A +X_1p& B_1+Y_1p  \\
 X_2p & Y_2p \\
  \end{array}
  \right)
: \left(\begin{array}{cc}
                                  X_1& Y_1  \\
                                  X_2  & Y_2 \\
                                \end{array}
                              \right)\in \mathcal{C}_{H}
   \right\}.$}$$
By (\ref{5HF713SFBlk3}), $\mbox{\small$\left(\begin{array}{cc}
A +X_1p& B_1+Y_1p  \\
 X_2p & Y_2p \\
  \end{array}
  \right)  \sim \left(\begin{array}{cc}
X_1p& B_1+Y_1p  \\
 0 & 0 \\
  \end{array}
  \right)$}$,  hence $\small\rho\left(
                                  \begin{array}{cc}
                                    A & 0 \\
                                    X_2p  & Y_2p
                                  \end{array}
                                \right)=\rho\left(
                                  \begin{array}{cc}
                                    A & 0 \\
                                    0  & Y_2p
                                  \end{array}
                                \right)\leq d-1$, which implies that $Y_2p=0$. By (\ref{5HF713SFBlk3}) again,
$\mbox{\small$\left(\begin{array}{cc}
A +X_1p& B_1+Y_1p  \\
 X_2p & 0 \\
  \end{array}
  \right) \sim \left(\begin{array}{cc}
X_1p & B_1+(0,e_i)+Y_1p  \\
 0 & 0 \\
  \end{array}
  \right)$}$, and hence
$\small\rho\left(
                                  \begin{array}{ccc}
                                    A & 0 & -e_i \\
                                    X_2p  & 0&0
                                  \end{array}
                                \right)=\rho\left(
                                  \begin{array}{ccc}
                                    A & 0 & -e_i\\
                                    0 &0 & X_2A^{-1}e_ip
                                  \end{array}
                                \right)\leq d-1$, $i=1,\ldots, d-1$. It follows that $X_2A^{-1}e_i=0$, $i=1,\ldots, d-1$, thus  $X_2=0$.
Therefore, we get that
\begin{equation}\label{uty7407dfvba523d}
\mbox{$\pi^{-1}\left(\begin{array}{ccc}
A & B_1  \\
 0 & 0  \\
  \end{array}
  \right)=
\left\{\left(\begin{array}{ccc}
A +Xp& B_1+Yp  \\
 0 & 0  \\
  \end{array}
  \right)
: (X, Y)\in\mathbb{Z}_{p^{s-1}}^{(d-1)\times n}\right\}$,}
\end{equation}
for all $A\in GL_{d-1}(\mathbb{Z}_{p})$ and $B_1\in\mathbb{Z}_p^{(d-1)\times(n-d+1)}$.

Finally, let $B\in \mathbb{Z}_p^{(d-1)\times(d-1)}$ and $B_1\in\mathbb{Z}_p^{(d-1)\times(n-d+1)}$. Write $L=(B,B_1)$.
 By (\ref{maxicliq5436uvbvl088}), there is a maximum  clique $\mathcal{C}_{L}$
of $\Gamma_d(\mathbb{Z}_{p^{s-1}}^{m\times n})$  such that
$$\mbox{$\pi^{-1}\left(\begin{array}{cc}
B & B_1 \\
0& 0  \\
  \end{array}
  \right)=
\left\{\left(\begin{array}{cc}
B +X_1p& B_1+Y_1p  \\
 X_2p & Y_2p \\
  \end{array}
  \right)
: \left(
                                \begin{array}{cc}
                                  X_1& Y_1  \\
                                  X_2  & Y_2 \\
                                \end{array}
                              \right)\in \mathcal{C}_{L}
   \right\}.$}$$
By Lemma \ref{rankinvertible935tet}, there is $A\in GL_{d-1}(\mathbb{Z}_{p})$ such that $B-A$ is invertible.
Applying (\ref{uty7407dfvba523d}),
$$\mbox{$\left(\begin{array}{cc}
B+X_1p& B_1+Y_1p  \\
 X_2p & Y_2p \\
  \end{array}
  \right)  \sim \left(\begin{array}{cc}
A +X_1p& B_1+Y_1p  \\
 0 & 0 \\
  \end{array}
  \right)$},$$
  and hence $\small\rho\left(
                                  \begin{array}{cc}
                                    B-A & 0 \\
                                    X_2p  & Y_2p
                                  \end{array}
                                \right)=\rho\left(
                                  \begin{array}{cc}
                                    B-A & 0 \\
                                    0  & Y_2p
                                  \end{array}
                                \right)\leq d-1$. Thus  $Y_2p=0$.
By (\ref{uty7407dfvba523d}) again, we have
$$\mbox{$\left(\begin{array}{cc}
B+X_1p& B_1+Y_1p  \\
 X_2p & 0 \\
  \end{array}
  \right) \sim \left(\begin{array}{cc}
A+X_1p & B_1+(0,e_i)+Y_1p  \\
 0 & 0 \\
  \end{array}
  \right)$}, \ \ i=1, \ldots, d-1, $$
  and hence
$$\rho\left(
                                  \begin{array}{ccc}
                                    B-A & 0 & -e_i \\
                                    X_2p  & 0&0
                                  \end{array}
                                \right)=\rho\left(
                                  \begin{array}{ccc}
                                    B-A & 0 & -e_i\\
                                    0 &0 & X_2(B-A)^{-1}e_ip
                                  \end{array}
                                \right)\leq d-1,  \ \ i=1, \ldots, d-1.$$
 Consequently, $X_2(B-A)^{-1}e_ip=0$, $i=1, \ldots, d-1$.
Thus it is clear that $X_2=0$. Therefore, we obtain that
\begin{equation}\label{6gfdre1vzkk3d}
\mbox{$\pi^{-1}\left(\begin{array}{ccc}
B & B_1  \\
 0 & 0  \\
  \end{array}
  \right)=
\left\{\left(\begin{array}{ccc}
B +Xp& B_1+Yp  \\
 0 & 0  \\
  \end{array}
  \right)
: (X, Y)\in\mathbb{Z}_{p^{s-1}}^{(d-1)\times n}\right\}$,}
\end{equation}
for all $B\in \mathbb{Z}_p^{(d-1)\times(d-1)}$ and $B_1\in\mathbb{Z}_p^{(d-1)\times(n-d+1)}$. Then
we have proved that
$$ {\cal M}=\pi^{-1}(\pi({\cal M}))=\left(
\begin{array}{c}
\mathbb{Z}_{p^s}^{(d-1)\times n} \\
0 \\
\end{array}\right).$$

{\bf Case 2.} \  $\pi({\cal M})$ is of the form (\ref{675fd425klhlcvx07}) with $n=m$.

Let $^t{\cal M}=\{ \,^tX: X\in {\cal M}\}$. Then $^t{\cal M}$ is a maximum clique of $\Gamma_d(\mathbb{Z}_{p^s}^{m\times m})$.
By (\ref{pipipipi004}), $\pi(^t{\cal M})=\,^t(\pi({\cal M}))$ is of the form (\ref{ghf435nvz9078rf}).
By Case 1, we have that $^t{\cal M}=\left(
\begin{array}{c}
\mathbb{Z}_{p^s}^{(d-1)\times m} \\
0 \\
\end{array}\right)$. Thus ${\cal M}=\left(\mathbb{Z}_{p^s}^{(m)\times (d-1)}, 0  \right)$.
$\qed$

\subsection{Cores of $\Gamma_d$ and the complement of $\Gamma_d$}

\ \ \ \ \
 A graph $G$ is a {\em core} \cite{Godsil} if every endomorphism of $G$ is an automorphism.
A subgraph $H$ of a graph $G$ is a {\em core of $G$} \cite{Godsil} if it is a core and there exists some homomorphism from $G$ to $H$.
Every graph $G$ has a core, which is an induced subgraph and is unique up to isomorphism \cite[Lemma 6.2.2]{Godsil}.

\begin{lem}\label{te65yer348nnv4}{\em (see  \cite[Lemma 2.5.9]{Roberson})} \ Let $G$ be a graph. Then the core of $G$ is the complete graph  $K_r$
if and only if $\omega(G)=r=\chi(G)$.
\end{lem}

\begin{cor}\label{main5350gdj67}
Let $n\geq m$. Then the core of $\Gamma_d(\mathbb{Z}_{p^s}^{m\times n})$ is a maximum clique of  $\Gamma_d(\mathbb{Z}_{p^s}^{m\times n})$,
and $\Gamma_d(\mathbb{Z}_{p^s}^{m\times n})$ is not a core.
\end{cor}
\proof
Let $\Gamma=\Gamma_d(\mathbb{Z}_{p^s}^{m\times n})$. By (\ref{gf353tgdcvadkhl54}), we have $\chi(\Gamma)=\omega(\Gamma)$.
By Lemma \ref{te65yer348nnv4},  the core of $\Gamma$ is a maximum clique, and hence $\Gamma_d(\mathbb{Z}_{p^s}^{m\times n})$ is not a core.
$\qed$

Let $\overline{G}$ denote the {\em complement} of a graph $G$. Then $\alpha(\overline{G})=\omega(G)$,
$\omega(\overline{G})=\alpha(G)$,  and ${\rm Aut}(\overline{G})={\rm Aut}(G)$. Moreover,  $\chi(\overline{G})\geq \omega(\overline{G})=\alpha(G)$.

\begin{thm}\label{complement-01}
Let $\Gamma=\Gamma_d(\mathbb{Z}_{p^s}^{m\times n})$ where $n\geq m$. Then  $\chi(\overline{\Gamma})=\alpha\left(\Gamma \right)=\omega(\overline{\Gamma})$.
Moreover, the core of $\overline{\Gamma}$ is a maximum clique of $\overline{\Gamma}$, and $\overline{\Gamma}$ is not a core.
\end{thm}
\proof
Let $\mathcal{S}=\left\{S_1, S_2, \ldots, S_{\alpha}\right\}$ be a largest independent set
of $\Gamma$, where $\alpha:=\alpha(\Gamma)=p^{sn(m-d+1)}$. Let
$$\mathcal{M}_1=\left\{\left(
                         \begin{array}{c}
                           X \\
                           0 \\
                         \end{array}
                       \right): X\in\mathbb{Z}_{p^s}^{(d-1)\times n} \right\}\subset \mathbb{Z}_{p^s}^{m\times n}.$$
By (\ref{biline13nvhiu6}), $\mathcal{M}_1$ is a maximum clique of  $\Gamma$.
Put $\mathcal{C}_i=\mathcal{M}_1+S_i$, $i=1,\ldots,\alpha$. Then $\mathcal{C}_1, \ldots, \mathcal{C}_{\alpha}$ are $\alpha$  maximum cliques of $\Gamma$,
and $\mathcal{C}_i\cap \mathcal{C}_j=\emptyset$ for all $i\neq j$. By Theorem \ref{bilinearformsyr56bch8},  we get $\left|V(\Gamma)\right|=p^{smn}=\alpha \cdot\omega(\Gamma)$.
Thus,
$V(\Gamma)$ has a  partition into $\alpha$  maximum cliques: $V(\Gamma)=\mathcal{C}_1\cup \mathcal{C}_2\cup \cdots \cup \mathcal{C}_{\alpha}$.
Since  $\mathcal{C}_1, \ldots, \mathcal{C}_{\alpha}$ are $\alpha$ largest independent sets of $\overline{\Gamma}$,
$V(\overline{\Gamma})$ has a  partition into $\alpha$  largest independent sets:
$V(\overline{\Gamma})=\mathcal{C}_1\cup \mathcal{C}_2\cup \cdots \cup \mathcal{C}_{\alpha}$. Hence $\chi(\overline{\Gamma})\leq \alpha$.
By $\chi(\overline{\Gamma})\geq \alpha=\omega(\overline{\Gamma})$, we obtain $\chi(\overline{\Gamma})=\alpha\left(\Gamma \right)=\omega(\overline{\Gamma})$.
By Lemma \ref{te65yer348nnv4}, we have similarly that  the core of $\overline{\Gamma}$ is a maximum clique,  and hence $\overline{\Gamma}$ is not a core.
$\qed$




\begin{thebibliography}{36}
\addtolength{\baselineskip}{8mm}
\footnotesize
\addtolength{\itemsep}{-0.8 em}


\bibitem{Bini} G. Bini, F. Flamini, Finite Commutative Rings and Their Applications, Springer Science+Business Media, LLC,  New York, 2002.


\bibitem{brown} W.C. Brown, Matrices over Commutative Rings, Marcel Dekker, Inc.,  New York, 1993.


\bibitem{Brouwera2} A.E. Brouwer, A.M. Cohen, A. Neumaier, Distance-Regular Graphs, Springer-Verlag, Berlin, Heidelberg, New York, 1989.


\bibitem{Chartrand} G. Chartrand, P. Zhang, Chromatic Graph Theory, Taylor \& Francis Group, Boca Raton, London,  2009.


\bibitem{Freeideal} P.M. Cohn, Free Ideal Rings and Localization in General Rings, Cambridge University Press, Cambridge, 2006.


\bibitem{cohn2} P.M. Cohn, Free Ring and Their Relations, Second edition, Academic Press, London, 1985.


\bibitem{Cossidente} A. Cossidente, G. Marino, F. Pavese, Non-linear maximum rank distance codes, Des. Codes Cryptogr.  79 (2016) 597-609.


\bibitem{Delsarte} P. Delsarte, Bilinear forms over a finite field, with applications to coding theory, J. Combin. Theory, Ser. A 25 (1978) 226-241.




\bibitem{Dougherty2} S.T. Dougherty,  K. Shiromoto, MDR Codes Over $\mathbb{Z}_k$, IEEE Trans. Inform. Theory 46 (2000) 265-269.

\bibitem{Oued} M. El Oued, On MDR codes over a finite ring, Int. J. Inform. Coding Theory   3 (2015) 107-119.







\bibitem{Gabidulin} E.M. Gabidulin, Theory of codes with maximum rank distance, Probl. Inf. Transm. 21(1985) 1-12.

\bibitem{Gadouleau} M. Gadouleau, Z. Yan, Constant-rank codes, IEEE Trans. Inform. Theory  56 (2008) 876-880.

\bibitem{Godsil3} C. Godsil,   Interesting graphs and their colourings, Unpublished notes, 2004.

\bibitem{Godsil} C. Godsil,  G.  Royle,  Algebraic Graph Theory,  Springer, New York, 2001.




\bibitem{Horlemann} A.-L. Horlemann-Trautmann, New Criteria for MRD and Gabidulin codes and some rank-metric code constructions, arXiv: 1507.08641v3 [cs IT], 2016.



\bibitem{Huang-2014} L.-P. Huang, Z.-J. Huang, C.-K. Li and N.-S. Sze,  Graphs associated with matrices over finite fields and their endomorphisms,
Linear Algebra  Appl.  447 (2014) 2-25.


\bibitem{Huang-2017} L.-P. Huang, H.D. Su, G.H. Tang, J.-B. Wang, Bilinear forms graphs over residue class rings, Linear Algebra Appl. 523 (2017) 13-32.


\bibitem{T.Huang} T. Huang, An analogue of the  Erd${\rm\ddot{o}}$s-Ko-Rado theorem for the distance-regular graphs of bilinear forms, Discrete Math. 64 (1987) 191-198.


\bibitem{Lee} H. Lee, Y. Lee, Construction of self-dual codes over finite rings $\mathbb{Z}_{p^m}$, J. Combin. Theory, Ser. A  115 (2008) 407-422.





\bibitem{mcdonald1} B.R. McDonald, Finite Rings with Identity, Marcel Dekker, Inc., New York, 1974.



\bibitem{Newman} M. Newman, Integral Matrices, Academic Press, New York, London, 1972.

\bibitem{Neri} A. Neri, A.-L. Horlemann-Trautmann, T. Randrianarisoa, J. Rosenthal,  On the genericity of maximum rank distance and
Gabidulin codes, Des. Codes Cryptogr. DOI 10.1007/s10623-017-0354-4,  Published online: 08 April 2017.


\bibitem{Roberson} D.E. Roberson, Variations on a Theme: Graph Homomorphisms, Ph.D. dissertation, University of Waterloo, 2013.


\bibitem{Roth} R.M. Roth, Maximum-rank array codes and their application to crisscross error correction, IEEE Trans. Inf. Theory 37(1991) 328-336.


\bibitem{Roth2} R.M. Roth, Introduction to Coding Theorey (Reprinted with corrections), Cambridge University Press, Cambridge, 2007.






\bibitem{Tanaka} H. Tanaka, Classification of subsets with minimal width and dual width in Grassmann, bilinear forms and dual polar graphs, J. Combin. Theory, Ser. A 113 (2006) 903-910.


\bibitem{GaloisRing} Z.-X. Wan,  Lectures on Finite Fields and Galois Rings, World Scientific, New Jersey, London, Singapore,  Hong Kong, 2003.





\bibitem{Y.X.Wang} Y.-X. Wang, Y.-J. Huo, C.-L. Ma,  Association Schemes of Matrices,  Science Press, Beijing, 2011.




\end{thebibliography}
\end{document}